\documentclass[12pt]{amsart}
\usepackage{amssymb, amsthm, amsmath}
\usepackage[english]{babel}
\usepackage{amsfonts}
\DeclareMathOperator{\id}{id}

\newcommand{\ra}{\rightarrow}

\newcommand{\ot}{\otimes}
\newcommand{\mtc}{\mathcal}
\newcommand{\lam}{\lambda}

\newcommand{\Lam}{\Lambda}

\newcommand{\al}{\alpha}
\newcommand{\eps}{\epsilon}
\newcommand{\bn}{\begin}
\newcommand{\en}{\end}

\newcommand{\mc}{\mathcal}

\newcommand{\nhs}{normal Hopf subalgebra}

\numberwithin{equation}{section}
\newtheorem{lem}[equation]{Lemma}
\newtheorem{thm}[equation]{Theorem}
\newtheorem{prop}[equation]{Proposition}
\newtheorem{defn}[equation]{Definition}
\newtheorem{cor}[equation]{Corollary}
\newtheorem{rem}[equation]{Remark}
\newcommand{\ch}{\chi}
\newcommand{\mtr}{\mathrm}

\numberwithin{equation}{section}

\newcommand{\ncm}{\newcommand}
\ncm{\np}{\newpage}
\ncm{\ebl}{\end{thebibliography}}
\ncm{\bbl}{\begin{thebibliography}}
\ncm{\chd}{_{ _{\ch}}}
\ncm{\ald}{_{ _{\al}}}
\ncm{\cP}{\mathcal{P}}
\ncm{\ei}{e_i}
\ncm{\eij}{e_{i,\;j}}
\ncm{\bt}{\begin{thm}}

\ncm{\bdef}{\begin{defn}
}

\ncm{\edf}{\end{defn}}
\ncm{\et}{\end{thm}}
\ncm{\bc}{\begin{cor}}
\ncm{\bl}{\begin{lem}}
\ncm{\el}{\end{lem}}
\ncm{\bpf}{\begin{proof}}
\ncm{\epf}{\end{proof}}
\ncm{\ec}{\end{cor}}
\ncm{\er}{\end{rem}}
\ncm{\br}{\begin{rem}}
\ncm{\bp}{\begin{prop}}
\ncm{\ep}{\end{prop}}
\ncm{\bd}{\begin{document}}
\ncm{\ed}{\end{document}}
\ncm{\beq}{\begin{equation}}
\ncm{\beqn}{\begin{equation*}}
\ncm{\eeq}{\end{equation}}
\ncm{\eeqn}{\end{equation*}}
\ncm{\bea}{\begin{eqnarray}}
\ncm{\eea}{\end{eqnarray}}
\ncm{\beanon}{\begin{eqnarray*}}
\ncm{\eeanon}{\end{eqnarray*}}

\title[Hopf algebras]
{On complements and the factorization problem for Hopf algebras}
\author{Sebastian Burciu}
\address{Inst.\ of Math.\ ``Simion Stoilow" of the Romanian Academy
P.O. Box 1-764, RO-014700, Bucharest, Romania}\address{\textnormal{and}} \address{University of Bucharest, Faculty of Mathematics and Computer Science, 14 Academiei St., Bucharest, Romania }
\email{sebastian.burciu@imar.ro}
\date{\today}

\begin{document}
\subjclass[2000]{Primary 16W30}
\keywords{Normal Hopf subalgebras; Semi-direct product; Internal tensor product}

\begin{abstract}
Two new results concerning complements in a semisimple Hopf algebra are proved. They extend some well known results from group theory. The uniqueness of Krull Schmidt Remak type decomposition is proved for semisimple completely reducible Hopf algebras.
\end{abstract}

\maketitle
\section{Introduction}

In recent years many concepts and results from the theory of finite groups
have been generalized or adapted to the context of (semisimple) Hopf algebras. Simultaneously, many results on the classification of semisimple Hopf algebras have also appeared. However some results, such as Kaplansky's conjecture about the  dimensions of the irreducible modules, still remain an open problem. An important class of mysterious Hopf algebras is that of simple Hopf algebras, i.e not having normal Hopf subalgebras. It was proven for example that if $H$ is a non-trivial semisimple Hopf algebra of dimension $< 60$, then $H$ is not simple (see \cite{ama-sonia}). Thus the smallest example of a simple Hopf algebra is in dimension $60$.

Let $H$ be a Hopf algebra and $A$ be a normal Hopf subalgebra of $H$. We say that $A$ has a {\it complement} in $H$ if and only if there is a Hopf subalgebra $L$ of $H$ such that $H \cong A\otimes  L$ as Hopf algebras. In this situation we also say that $A$ is a {\it direct factor} of $H$. Also we say that $(A, L)$ is a {\it normal factorization} of $H$.

Recall that $(A,L)$ is a factorization of $H$ if the multiplication map $$m:A\otimes L \ra H,\;\;a\otimes l \mapsto al $$ is a vector space isomorphism.

In this paper we prove two results on factorization of Hopf algebras. They are analogous results to some results from group theory. The first result shows that a given factorization with both terms normal Hopf subalgebras is a normal factorization. The second results concerns the case of factorizations $(A,\;L)$ where just one of the terms, say $A$, ia a normal Hopf subalgebra. In this situation we show that $H \cong A \# L$ as Hopf algebras.

Next we apply the previous results to the case of semisimple Hopf algebras. In this situation we can use Theorem \ref{dimform} proven in \cite{Bd} to obtain new results for semisimple Hopf algebras. For example, using these results, similar arguments to those from group theory can be applied to obtain a characterization of completely reducible Hopf algebras. Parallel to the theory of groups we call a Hopf algebra $H$ {\it completely reducible} if it can be written as an internal direct tensor product:
\begin{equation}
  H=\otimes_{ \lam \in \Lam}K_{\lam},
\end{equation}
of the simple semisimple Hopf algebras $K_{\lam}$.
Definition of internal tensor product is given in Section 4. For $\Lam$ a finite set this coincides with the usual tensor product, see Propostion \ref{int}.

A Krull-Schmidt-Remak type uniqueness for such a decomposition is proven in the finite dimensional case.

The paper is organized as follows. In the second section we consider factorizations of Hopf algebras. If one factor of the factorization is normal then we show in Theorem \ref{main1} that one obtains a semi-direct product of Hopf algebras. If both factors are normal we show in Theorem \ref{main} that in this case one obtains the usual tensor product of Hopf algebras. In the proof of both theorems we use the fact a factorization is the same as a bicrossed product Hopf algebra \cite{Majid}.

The third section applies the results of the previous section to semisimple Hopf algebras. The key point of this section is to use Theoerem \ref{dimform} on the dimension of $AL$ inside $H$ in terms of the dimensions of $A$, $L$ and of the intersection $A \cap L$.

The next section considers completely reducible Hopf algebras and shows the uniqueness of a  decomposition of Krull-Schmidt-Remak type  for finite dimensional completely reducible Hopf algebras.

In the last section we define similar to the group theory the socle of a semisimple Hopf algebra and study the lattice of normal Hopf subalgebras of a semisimple Hopf algebra.

We use the standard notations for Hopf algebras that can be found for example in \cite{Montg}. We work over an algebraically closed field $k$ of characteristic zero.
%%%%%%%%%%%%
\section{Factorization of Hopf algebras}

In this section we prove two results on factorization of arbitrary Hopf algebras.
\subsection{Settings}
Let $H$ be a Hopf algebra over a field $k$. Recall that a Hopf subalgebra $A \subseteq H$ is called \emph{normal} if $h_1 A S(h_2) \subseteq A$ and $S(h_1)Ah_2\subseteq A$, for all $h \in H$. If $H$ does not contain proper normal Hopf subalgebras then it is called \emph{simple}. Recall from \cite{Montg} that $A$ is normal in $H$ if and only if $HA^+=A^+H$ and in this situation $H//A: =H/HA^+$ is a quotient Hopf algebra of $H$. Here $A^+:=A \cap \mtr{ker}\;\eps$.
 %
 %This is equivalent to $Sh_1Ah_2 \subset A$ for all $h \in H$ (see Lemma 1.1 of \cite{mascoid}).
 %
\bn{lem}\label{arg}
Suppose that $A$ and $L$ are Hopf subalgebras of $H$ with $A$ a normal Hopf subalgebra of $H$. Then $AL=LA$ is a Hopf subalgebra of $H$.
\end{lem}
\bpf
Indeed, note that $la=(l_1aSl_2)l_3 \in AL$ for all $l \in L$ and $a \in A$. Thus $LA \subset AL$. Similarly $al=l_1(Sl_2al_3) \in LA$ which shows that $AL\subset LA$.
\epf

\subsection{Bicrossed product} The bicrossed product of Hopf algebras was
introduced by Majid in \cite{Majid}, initially under the name of double cross product.

Recall that a matched pair of bialgebras is a system $(A,L,\lhd,\rhd)$, where $A$ and $L$ are bialgebras, $\lhd : L \otimes A \ra L$ and $\rhd : L\otimes A \ra A$ are coalgebra maps such that $(A,\rhd)$ is a left $L$-module coalgebra, $(L,\lhd)$ is a right $A$-module coalgebra and the following compatibility conditions hold:

\begin{equation*}
1_L \lhd a = \eps_A(a)1_L,\;\; l \rhd 1_A = \eps_L(l)1_A
\end{equation*}

\begin{equation*}
l \rhd (ab) = (l_1 \rhd a_1)((l_2 \lhd a_2) \rhd  b)
\end{equation*}

\begin{equation*}
    (lh) \lhd a = ((l \lhd h_1) \rhd a_1)(h_2 \rhd a_2)
\end{equation*}

\begin{equation*}
    l_1 \lhd a_1 \otimes l_2 \rhd a_2 = l_2 \lhd a_2 \otimes l_1 \rhd a_1
\end{equation*}
for all $a, b \in A,\; l, h \in L$. If $(A,L,\lhd,\rhd)$ is a matched pair of bialgebras then the bicrossed product $A \bowtie L$ of $A$ with $L$ is the vector space $A \otimes L$ endowed with the multiplication given by
\begin{equation}\label{bicross}
    (a \bowtie l) · (b \bowtie h) := a(l_1 \rhd b_1) \bowtie (l_2 \lhd b_2)h
\end{equation}
for all $a, b \in A,\; l, h \in L$, where $a \otimes l$ is denoted by $a \bowtie l$. Then $A \bowtie L$ is a bialgebra with the coalgebra structure given by the tensor product of coalgebras. Moreover, if $A$ and $L$ are Hopf algebras, then $A \bowtie L$ has an antipode given by the formula:
$$S(a \bowtie l) := (1_A \bowtie S_L(l)) · (S_A(a) \bowtie 1_L) $$
for all $a \in A$ and $l\in L$ [ see \cite{Majid}, Theorem 7.2.2].
\subsection{Factorization of Hopf algebras}
Let $H$ a Hopf algebra and $A$, $L$ be two Hopf subalgebras. We say that $(A,L)$ is a {\it factorization} of $H$ if the multiplication map $$m: A \otimes L \ra H,\;\; a \ot l \ra al $$ is bijective. Note that in this case $A \cap L= k$.

Then by Theorem 2.7.3 from \cite{Majid} it follows that $H$ is a bicrossed product Hopf algebra of $A$ and $L$. The bicrossed product setting is recovered as following:

One defines the $k$-linear map
\begin{equation*}
    \mu: L \otimes A \ra A \otimes L, \;\;\mu(l \otimes a) := u^{-1}(la)
\end{equation*}
for all $l \in L$ and $a \in A$. Then the actions $\lhd,\; \rhd$ are given by the formulae:

\begin{equation}\label{rhd}
\rhd : L \otimes A \ra A, \;\; l \rhd a := ((\mtr{Id} \otimes \eps_L) \circ \mu)(l \otimes a)
\end{equation}

\begin{equation}\label{lhd}
\lhd : L \otimes A \ra L, \;\;h \lhd a := ((\eps_A \otimes \mtr{Id}) \circ \mu)(l \otimes a) .
\end{equation}

\subsection{Factorization with a normal factor}

In this subsection we prove that  for a factorization  $(A,\;L)$ of $H$ with a \nhs\; $A$ of $H$ one has that $H$  is a semi-direct product of $A$ and $L$ as $k$-algebras.

\bn{prop}\label{main1}
Suppose that  $(A,L,\lhd,\rhd)$ is a bicrossed product of $H$ with $A$ a normal Hopf subalgebra of $L$. Then $A \# L\cong H $ as Hopf algebras via the multiplication map $\phi$ given by $a \bowtie l \mapsto al$.
\end{prop}
% definition of semidirect product
\bpf
Let $L$ acting on $A$ by the adjoint action $l.a=l_1aS(l_2)$. Since $A$ is a \nhs\; of $H$ it follows that $\mu(l \ot a)= l_1aS(l_2)\otimes l_3$ for all $l \in L$ and $a \in A$. Then the action $\lhd$ becomes trivial since $l \lhd a := ((\eps_A \otimes \mtr{Id}) \circ \mu)(l \otimes a)=(\eps_A \otimes \mtr{Id}) (l_1aSl_2\otimes l_3)=\eps(a)l$.

Then formula \ref{bicross} for the algebra product shows that $H \# A\cong L $ as Hopf algebras via the map $\phi$.
\epf

\bn{cor}\label{qouti}
Suppose that $H$ is a factorization by $A$ and $L$. If $A$ is a normal Hopf subalgebra of $H$ then the map $\phi: L \ra H//A$ given by $ l \mapsto \bar{l}$ is an isomorphism of Hopf algebras.%\psi
\end{cor}

\bn{proof}
Clearly $\phi$ is a Hopf algebra map. Define the inverse map of $\phi$ by $$\psi: H//A \ra L,\;\;\psi(\overline{a \bowtie l})=\eps(a)l.$$ It is easy to see that the map $\psi$ is well defined. On the other hand clearly $\psi$ is the inverse of the map $ \phi$.
\end{proof}

\subsection{Factorizations with both factors normal}
In this subsection we prove that for a factorization $(A,\;L)$ of $H$ with both $A$ and $L$ \nhs s\; of $H$ one has that $H$ is isomorphic as Hopf algebras with the tensor product $A \otimes L$.

\bn{prop}\label{commute}
Suppose that $M$ and $N$ are normal Hopf subalgebras of a given Hopf algebra $H$ such that $M \cap N=k$. Then $mn=nm$ for all $m \in M$ and $n \in N$.
\end{prop}

\bn{proof}
Note that $m_1n_1S(m_2)S(n_2) \in M\cap N$. Thus $m_1n_1S(m_2)S(n_2)=\eps(m)\eps(n)1$. This implies that

\begin{eqnarray*}
mn &=& m_1n_1S(m_2)S(n_2)n_3m_3 \\
= \eps(n_1)\eps(m_1)n_2m_2 &=& nm
\end{eqnarray*}
\end{proof}
%{\bf this works with coideal subalgebras}

\bt\label{main}
Suppose that  $(A,L,\lhd,\rhd)$ is a bicrossed product of $H$ and that both $A$ and $L$ are normal Hopf subalgebras of  $H$. Then $H\cong A\ot L$ as Hopf algebras via the multiplication map $a \ot l \mapsto al$.
\et

\bpf
Since $A \cap L=k$ previous proposition implies that $al=la$ for all $a \in A$ and $l \in L$. Then clearly $\mu(l \ot a)=a\ot l$ and formulae \ref{rhd} and \ref{lhd} show that both actions $\rhd$ and $\lhd$ are trivial. In this situation fromula \ref{bicross} show that $H$ is the tensor product of $A$ and $L$.
\epf

\section{Factorization of semisimple Hopf algebras}

In this section we apply the results from the previous sections in the case that $H$ is a semisimple Hopf algebra.

The following theorem proven in \cite{Bd} is needed in the sequel.

\begin{thm} \label{dimform}
Let $H$ be a semisimple Hopf algebra and $A,\;L$ be two Hopf subalgebras of $H$. Then
\begin{equation}
    L//L\cap A\cong LA//A
\end{equation}
 as coalgebras.
\en{thm}

Recall that $LA//A:=LA/LA^+$ is just a coalgebra and not a Hopf algebra if $A$ is not normal Hopf subalgebra of $H$.
From this it follows
\begin{equation}\label{dim}
    |LA|=\frac{|L||A|}{|L\cap A|}
\end{equation}
Here  for a vector space $V$ we denoted by $|V|$ its dimension as vector space.

Next we want to prove the following theorem, a result similar to the well known result on normal complements for groups:

\bn{cor}\label{main1ss}
Suppose that $A$ and $L$ are semisimple Hopf subalgebras of a Hopf algebra $H$ such that $A\cap L=k$. Moreover suppose that $A$ is a normal Hopf subalgebra of $H$. Then $A \# L\cong AL$ as Hopf algebras via the multiplication map $\phi$ given by $a \ot l \mapsto al$.
\end{cor}

\bpf
By Lemma \ref{arg} one has that $AL$ is a Hopf algebra subalgebra of $H$.

Define the adjoint action of $L$ on $A$ by $l.a=l_1aS(l_2)$ and form the smash product $A \# L$ (see \cite{Kas}). It can be easily be checked that the map $\phi$ defined above is a surjective algebra map. Indeed, to check that $\phi$ is an algebra map one has that

\begin{eqnarray*}
\phi((a\# l)(a'\# l')) &=& \phi(a(l_1a'Sl_2)\# l_3l')=\\
  a(l_1a'Sl_2)l_3l'=(al)(a'l') &=& \phi(a \# l)\phi(a' \# l')
\end{eqnarray*}
Since $A$ and $L$ are semisimple it follows that $A L$ is also semisimple Hopf algebra. Indeed by Lemma 2.16 of \cite{mascoid} it follows that $\Lam_A$ is central in $H$. Then it can be checked that $\Lam_A\Lam_L$ is an idempotent integral of $AL$ where $\Lam_L$ is the idempotent integral of $L$. Thus $AL$ is a semisimple Hopf algebra and from Theorem \ref{dimform} applied to $H:=AL$ it follows that $AL$ and $A \# L$ have the same dimension. This implies that $\phi$ is also injective and therefore an isomorphism of Hopf algebras by Theorem \ref{main1}.
\epf

As in the group situation one has the following result:

\bt\label{mainss}
Suppose that both $A$ and $L$ are normal Hopf subalgebras of a semisimple Hopf algebra $H$ such that $H=AL$ and $A\cap L=k$. Then $H\cong A\ot L$ as Hopf algebras via the multiplication map $a \ot l \mapsto al$.
\et

\bpf
As before it follows that $AL=LA$. Proposition \ref{commute} implies that $al=la$ for all $a \in A$ and $l \in L$. Define $\phi: A \otimes L  \ra AL$ by $a \otimes l \mapsto al$. Note that $\phi$ is a surjective algebra map. Clearly $\phi$ is a coalgebra map. Since by formula \ref{dim} one has that $|AL|=|A||L|$ it follows that $\phi$ is an isomorphism of Hopf algebras.
\epf

%{\bf Same thing with coideal subalgebras!}
\subsection{Projection maps}
%for nonsemisimple case too??
The treatment of this subsection follows \cite{Rob}.

Suppose that $H \cong A \otimes L$ is an exact normal factorization of a semisimple Hopf algebra $H$. Define the following projection maps $\pi_1 :H \ra H$ given by $\pi_1(a \ot l)=\eps(l)a$ and $\pi_2 :H \ra H$ given by $\pi_2(a \ot l)=\eps(a)l$. Note that $\pi_1$ and $\pi_2$ are Hopf algebra maps.

An endomorphism of $H$ is called {\it normal} if it is a morphism of modules under the adjoint action of $H$ on itself. Thus $\phi:H \ra H$ is normal if and only if $\phi(h_1aSh_2)=h_1\phi(a)Sh_2$ for all $a, h \in H$.
% Do I need semisimple?

\bp With the above notations, for $i=1,\; 2$, the maps $\pi_i$ satisfy the following:

1. $\pi_i$ are Hopf algebra maps.

2. $\pi_i$ are normal endomorphisms (morphisms of adjoint modules).

3. $\pi_i^2=\pi_i$ under composition.

3. $\pi_2\star \pi_1=\pi_1\star \pi_2= \id_{H}$ under convolution.

4. $\pi_2 \circ \pi_1=\pi_1\circ \pi_2=\eps(\;)1$.
\ep

\bpf
Straightforward computations.
\epf

\bt
Let $H$ be a semisimple Hopf algebra. There is a bijection between exact normal factorizations of $H \cong A \ot L$ with $A,\;L$ semisimple Hopf algebras and the set of pairs of two orthogonal normal idempotents $\{\pi_1,\;\pi_2\}$ satisfying the above properties.
\et

\bpf
In the previous proposition we have shown how to associate a pair of normal endomorphisms with the above properties to any exact normal factorization of $H$.
Conversely, given $\pi_i$ as above let $A:=\pi_1(H)$ and $L:=\pi_2(H)$. Then $A$ and $L$ are normal Hopf subalgebras since $\pi_1$ and $\pi_2$ are normal endomorphisms. Note that $A\cap L=k$. Indeed if $x \in A \cap L$ then $x=\pi_1(h)=\pi_2(h')$ for some $h,h' \in H$. Thus $\pi_1(x)=\pi_1^2(h)=\pi_1(h)=x$ and on the other hand $\pi_1(x)=\pi_1\pi_2(h')=\eps(h')1$ which implies that $x$ is a scalar. Next $h=\pi_1\circ \pi_2(h)=\sum\pi_1(h_1)\pi_2(h_2) \in AL$ which shows that $H \subset AL$. Theorem \ref{mainss} implies that $H \cong A \otimes L$ is a normal factorization.
\epf

\section{Completely reducible Hopf algebras}

In this section we prove few results on semisimple completely reducible Hopf algebra.

For a set of Hopf subalgebras $\{K_{\lam}\}_{ \lam \in \Lam}$ of $H$ denote by $< K_{\lam}\;|\; \lam \in \Lam>$ the Hopf subalgebra of $H$ generated by all these Hopf subalgebras. It coincides with the subalgebra of $H$ generated by $\{K_{\lam}\}_{ \lam \in \Lam}$. In analogy with group theory we give the following definition of internal tensor product.

\bn{defn}
We say that $H$ is the internal tensor product of the family of Hopf subalgebras $\{K_{\lam}\}_{ \lam \in \Lam}$ and write

$$H =\ot_{\lam \in \Lam}K_{\lam}$$

if $H=< K_{\lam}\;|\; \lam \in \Lam>$ and $K_{\mu}\;\cap \;< K_{\lam}\;|\; \lam \neq \mu>=k$ for all $\mu \in \Lam$.
\end{defn}

Nest result shows that for semisimple Hopf algebras finite internal tensor products coincide with the usual tensor products.

\bn{prop}\label{int}
Suppose that $\Lam=\{\lam_1,\cdots ,\lam_s\}$ is a finite set and $$H =\ot_{\lam \in \Lam}K_{\lam}$$ is a finite internal tensor product of a semisimple Hopf algebra $H$. Then $$H \cong K_{\lam_1}\ot\cdots \ot K_{\lam_s}$$
as Hopf algebras.
\end{prop}

\bn{proof}
Since $K_{\lam_1}\cap <K_{\lam}\;|\; \lam \neq \lam_1>=k$ by Theorem \ref{mainss} it follows that $H\cong K_{\lam_1}\ot<K_{\lam}\;|\; \lam \neq \lam_1>$ as Hopf algebras. Inductively one can see that $<K_{\lam}\;|\; \lam \neq \lam_1>\cong K_{\lam_2}\ot\cdots \ot K_{\lam_s}$. Thus $H \cong K_{\lam_1}\ot\cdots \ot K_{\lam_s}$.
\end{proof}

\bn{rem}\end{rem}
It will be interesting to decide wether the internal tensor product as defined above coincides with the usual tensor product if $\Lam$ is infinite or if at least one of the Hopf subalgebras $K_{\lam}$ has infinite dimension. Certainly this is the case for groups (see \cite{Rob}).

\bn{defn} A Hopf algebra $H$ is called {\it completely reducible} if it can be written as a direct tensor product
$$
H=\otimes_{ \lam \in \Lam}K_{\lam}
$$
where $K_{\lam}$ are simple semisimple Hopf algebras. Such a decomposition will be called a {\it Krull-Schmidt-Remak type decomposition} of $H$ in analogy with the group situation.
\end{defn}

\bn{thm}\label{crd}
Suppose that $$
H=\otimes_{ \lam \in \Lam}K_{\lam}
$$ is a completely reducible semisimple Hopf algebra and $K$ is a normal Hopf subalgebra of $H$. Then there is $P \subset \Lam$ such that
$$
K=\otimes_{ \lam \in P}K_{\lam}.
$$
\end{thm}

\bn{proof} Since $H$ is semisimple it follows that $\Lam$ is a finite set.
One may suppose that $K \neq H$. Let $\mtc{S}$ be the following set:
\begin{equation*}
    \mtc{S}=\{ \Lam'\subset \Lam \;|\; <K, K_{\lam}\;| \lam \in \Lam'>\cong K \ot (\otimes_{ \lam \in \Lam'}K_{\lam} )\}
\end{equation*}
The set $\mtc{S}$ is not empty. Indeed since $H \neq K$ there is $\lam$ such that $K_{\lam}\nsubseteq K$. Since $K_{\lam}$ is simple it follows that $K\cap K_{\lam}=k$ and Theorem \ref{mainss} implies that $\{\lam\} \in \mtc{S}$. Since $\mtc{S}$ is finite it has a maximal element denoted by $M$. Let $L:=<K,\; K_{\lam}\;| \lam \in M>\cong K \ot (\otimes_{ \lam \in M}K_{\lam} )>$. If $\lam \in \Lam\setminus M$ then $L \cap K_{\lam} \
\neq k$ since $M\cup \{\lam\} \notin \mtc{S}$. Since $K_{\lam}$ is a simple Hopf algebra it follows that $K_{\lam} \subset L$ and thus $L=H$.

Put $P=\Lam \setminus M$. It follows that $K':=\otimes_{ \lam \in P}K_{\lam} \subset K$. Then applying Lemma \ref{qouti} the composition of the maps $K' \ra H//(\otimes_{ \lam \in M}K_{\lam}) \ra K$ given by $x \mapsto \bar{x} \mapsto x$ is the composition of two Hopf algebra isomorphisms and therefore it is a Hopf algebra isomorphism too.
\end{proof}

The previous theorem implies the following Corollary, that a finite Krull-Schmidt-Remak type decomposition is unique up to a permutation.

\bn{cor}
Suppose that $H=\otimes_{i=1}^rK_{i}$ and $H=\otimes_{j=1}^sL_{j}$ are two Krull-Schmidt-Remak  type decompositions of $H$. Then $r=s$ and there is a permutation $\sigma$ of $\{1,\cdots, r\}$ such that $K_i\cong L_{\sigma(i)}$.
\end{cor}

\bn{proof}
By previous Theorem $K_1\cong \otimes_{ j \in \Lam_1}L_j$ for some set $\Lam_1$. Since $K_1$ is simple it follows that $\Lam_{1}=\{\sigma(1)\}$ for some element $\sigma(1)\in \{1,\cdots, s\}$. Then one can continue arguing in the same way for all the other Hopf subalgebras $K_i$.
\end{proof}
Recall from introduction that a direct factor of $H$ is a \nhs \; $K$ of $H$ such that there is a normal Hopf subalgebra $L$ of $H$ with $(K,\;L)$ an exact normal factorization of $H$.

Next Theorem can be regarded as a converse of the previous Theorem. It generalizes Head's theorem from groups (see 3.3.13 of \cite{Rob}).

\bn{rem}\label{finite}\end{rem}
The lattice of normal Hopf subalgebras of a semisimple Hopf algebra is finite. This can be deduced from Theorem 2. 4 of \cite{Bker} for example. From this theorem to any normal Hopf subalgebra $K$ one can associate uniquely a central character $$t_K=\sum_{\{\ch \in \mtr{Irr}(H)\;|\;\mtr{ker}_H(\ch) \supset K\}}\ch(1)\ch.$$

Conversely recall the notion of kernel of a character introduced in \cite{Bker}. To any such central character $t$ one can associate a normal Hopf subalgebra of $H$ namely $K=\mtr{ker}\;t$. Since $\mtr{Irr}(H)$ is finite it follows that the lattice of normal Hopf subalgebras of $H$ is also finite.

\bt
Let $H$ be a semisimple Hopf algebra such that any normal Hopf subalgebra of it is contained in a direct factor of $H$. Then $H$ is completely reducible.
\et

\bpf
Let $L$ be the Hopf subalgebra generated by all simple normal Hopf subalgebras of $H$. Then using previous Remark it follows that $L$ is completely reducible. Thus if $L=H$ we are done. Suppose that $L \neq H$ and let $a \in H\setminus L$. There is a normal Hopf subalgebra $M$ which is maximal with the property that  $L \subset M$ and $a \notin M$. %Indeed let $\mtc{S}$ be the set of all normal Hopf subalgebras with this property. Then $L \in \mtc {S}$ which shows that $\mtc{S}$ is not empty. Then one can see that any ascending chain of elements of $\mtc{S}$ has a maximal element, namely the union of its elements. Thus by Zorn's lemma it follows that $\mtc{S}$ has a maximal element.

The hypothesis imply that $M$ is contained in a direct factor $K$ with $H \cong K \ot K'$. Since $K'$ is not trivial it follows that $a \in M \ot K'$ by maximality of $M$. Also by formula \ref{dim} note that $K \cap (M\ot K')=M$ so $a \notin K$. Maximality of $M$ implies that $K=M$ and $H \cong M \ot K'$. If $K''$ is a normal Hopf subalgebra of $K'$ then it is easy to see that $K''$ is normal in $H$. Then $M \ot K"$ lies in a direct factor of $H$. But it was already have been shown that $M$ cannot be contained in any direct factor of $H$ except $M$ itself. Thus $K'$ is simple and  therefore $K' \subset L$ which is a contradiction. Thus $L=H$ is a completely reducible Hopf algebra.
\epf

%%%%%%%%%%%%%%%%%%%%%%%%%%%%%%%%%%%%%%%%%%%%
\section{Socle of a semisimple Hopf algebra}
The socle of a group is the subgroup generated by all minimal normal subgroups. Because the product of normal subgroups is a subgroup, it follows that the socle of a group is the direct product of some of its minimal subgroups. Using Theorem \ref{main} a similar fact can be proven for semisimple Hopf algebras.

\bn{lem}
Let $H$ be a finite dimensional semisimple Hopf algebra. If $M$ and $N$ are minimal normal Hopf subalgebras of $H$ then $M$ and $N$ centralize each other and therefore $MN \cong M \otimes N$.
\end{lem}

\bn{proof}
Since $N\cap M$ is a normal Hopf subalgebra of $H$, by the minimality of $M$ and $N$ it follows that $M \cap N=k$. Then one can apply Proposition \ref{commute}.
\end{proof}

In analogy with group theory one can define the {\it  socle} of a Hopf algebra as the Hopf subalgebra generated by all minimal normal Hopf subalgebras.

\bn{thm}
The socle of a semisimple Hopf algebras is the tensor product of some of its minimal normal Hopf subalgebras.
\end{thm}

\bn{proof}
Let $S$ be the socle of $H$. Let us assume that $L_1, \cdots , L_s$ are all minimal normal Hopf subalgebras of $H$. Then $L_1 \cap L_2=k $ and by Theorem \ref{mainss} one has $L_1L_2 \cong L_1\ot L_2$.

As $L_3$ is a minimal normal Hopf subalgebra and $L_{1} L_{2}$ is a normal Hopf subalgebra, $L_3$ is either contained in $L_{1} L_{2}$ or intersects it trivially. If $L_{3} $ is contained in $L_{1} L_{2}$ then skip it and consider the next minimal normal Hopf subalgebra $L_4$. Otherwise note that $L_1L_2L_3 \cong L_1 \ot L_2 \ot L_3$. Continuing this procedure we obtain a sequence $L_{i_1}\cdots L_{i_s}$ of minimal normal Hopf subalgebras where $S=L_{i_1}\cdots L_{i_s}$ and $L_{i_j}\cap (L_{i_1}\cdots L_{i_{j-1}})=k$.
It follows that $S \cong L_{i_1}\ot\cdots \ot L_{i_s}$.
\end{proof}

\subsection{On the lattice of normal Hopf subalgebras of a semisimple Hopf algebra}\label{normalfacts}
Let $H$ be a semisimple Hopf algebra. In this section we will show that there is an anti - isomorphism of the lattices of normal Hopf subalgebras of $H$ and $H^*$.

Suppose that $K$ is a normal Hopf subalgebra of $H$ and let $L=H//K$ be the quotient Hopf algebra of $H$ via $\pi:H\ra L$. Then $\pi^* :L^* \ra H^*$ is an injective Hopf algebra map. It follows that $\pi^*(L^*)$ is a normal  Hopf subalgebra of $H^*$. Indeed, it can be checked that \begin{equation}\label{forml}
\pi^*(L^*)=\{f\in H^*|f(ha)=f(h)\eps(a)\;\; ; h\in H,\;a\in K\}.
\end{equation}

Using this description it is easy to see that $g_1fS(g_2)\in L^*$ for all $g \in H^*$ and $f\in L^*$. Moreover ${(H^*//L^*)}^*\cong K$ since ${(H^*//L^*)}^*=\{a\in H^{**}=H\;|\; fg(a)=f(1)g(a) \;\; ; f\in H^* ,\; g \in L^*\}$.

{\bf }
\begin{thm} \label{dualcomp}Let $K$ and $L$ be normal Hopf subalgebras of a semisimple Hopf algebra $H$. Then
\bn{equation}
(H//(L K))^*=(H//L)^*\cap (H//K)^*
\end{equation}
and
\bn{equation}
(H//(L\cap K))^*=(H//L)^*(H//K)^*
\end{equation}

\end{thm}

\begin{proof}
The first equality follows directly using the characterization from Equation \ref{forml}. Indeed suppose that $f \in (H//(L K))^*$. Then by Equation \ref{forml} $f$ satisfies $f(ha)=f(h)\eps(a)$ for all $h\in H,\;a\in LK$. In particular $f(ha)=f(h)\eps(a)$ for all $h\in H,\;a\in K$ and  also $f(ha)=f(h)\eps(a)$ for all $h\in H,\;a\in L$. This shows that $f \in (H//L)^*\cap (H//K)^*$ and therefore
\bn{equation}
(H//(L K))^*\subset (H//L)^*\cap (H//K)^*
\end{equation}
The other inclusion is verified similarly.

The inclusion
\bn{equation*}
(H//L)^*(H//K)^*\subset(H//(L\cap K))^*
\end{equation*}
of the second item can also be checked using Equation \ref{forml}. Then it is enough to show that both terms of the second item have the same dimension, i.e:

\bn{equation}
|(H//(L\cap K))^*|=|(H//L)^*(H//K)^*|.
\end{equation}

Indeed, using previous item and formula \ref{dim} one has that

\begin{eqnarray*}
  |(H//L)^*(H//K)^*| &=& \frac{|(H//L)^*||(H//K)^*|}{|(H//L)^*\cap (H//K)^*|} \\
 &=& \frac{|(H//L)^*||(H//K)^*|}{|(H//LK)^*|}\\
 &=& \frac{|H|}{|L|}\frac{|H|}{|K|}/{\frac{|H|}{|LK|}} \\ &=& |H|/ \frac{|L||K|}{|LK|}\\
 &=& \frac{|H|}{|L\cap K|}\\ &=& |(H//(L\cap K))^*|
\end{eqnarray*}

 Thus
\bn{equation*}
(H//(L\cap K))^*=(H//L)^*(H//K)^*.
\end{equation*}
\end{proof}

For a Hopf algebra $H$ denote by $\mtc{L}(H)$ the lattice of normal Hopf subalgebras of $H$ with the operations of intersection and product. The previous theorem implies the following Corollary:

\bn{cor}
Let $H$ be a semisimple Hopf algebra over an algebraically closed field $k$. There is an anti-isomorphism $ K \mapsto (H//K)^*$ between $\mc{L}(H)$ and $\mtc{L}(H^*)$. Its inverse is given by $K \mapsto (H^*//K)^*$.
\end{cor}

\bn{rem}\end{rem}
Description of maximal normal Hopf subalgebras as kernels $\mtr{ker}\;\ch$ of central characters $\ch$ in $H^*$ was given in \cite{Bker}. The previous Corollary implies that the minimal normal Hopf subalgebras of $H^*$ are of the type $(H//\mtr{ker}\;\ch)^*$.

\subsection*{Aknowledgmets}

The author thanks A. Agore and G. Militaru for useful commments on a previous version of this manuscript.

This work was supported by the strategic grant \\ POSDRU/89/1.5/S/58852, Project "Postdoctoral programme for training scientific researchers" cofinanced by the European Social Found within the Sectorial Operational Program Human Resources Development 2007 - 2013.

\end{document}